\documentclass[a4paper]{article}

\usepackage{amsmath,amssymb,amsthm,algorithm}
\usepackage[noend]{algpseudocode}
\usepackage{hyperref}
\usepackage{tikz}

\newcommand{\Z}{\mathbb{Z}}
\newcommand{\N}{\mathbb{N}}

\newcommand{\realm}{\mathfrak{D}}

\newcommand{\lang}{\mathcal{L}}
\newcommand{\glimset}{\tilde \omega}
\newcommand{\limset}{\Omega}

\newcommand{\asymset}{\omega}

\newcommand{\set}[2]{\left\{ #1 \;\middle|\; #2 \right\}}
\newcommand{\uniform}[1]{#1^\Z}
\newcommand{\xt}[1]{\mathtt{#1}}

\theoremstyle{plain}
\newtheorem{lemma}{Lemma}

\newtheorem{theorem}[lemma]{Theorem}

\newtheorem{example}[lemma]{Example}

\title{Generically Nilpotent Cellular Automata}
\author{Ilkka T\"orm\"a\footnote{Research supported by Magnus Ehrnrooth Foundation.} \\ Department of Mathematics and Statistics \\ University of Turku, Finland \\ \texttt{iatorm@utu.fi}}

\begin{document}

\maketitle

\begin{abstract}
  We study the generic limit sets of one-dimensional cellular automata, which intuitively capture their asymptotic dynamics while discarding transient phenomena.
  As our main results, we characterize the automata whose generic limit set is a singleton, and show that this class is $\Sigma^0_2$-complete.
  We also prove that given a CA whose generic limit set is guaranteed to be a singleton, the sole configuration it contains cannot be algorithmically determined.
\end{abstract}

\section{Introduction}

A one-dimensional cellular automaton (CA for short) is a dynamical system whose phase space consists of bi-infinite sequences $x \in A^\Z$, where $A$ is a finite alphabet.
Every CA $f : A^\Z \to A^\Z$ is defined by a local rule that is applied to every coordinate simultaneously.
We consider cellular automata as topological dynamical systems.
In general, they are neither injective nor surjective, and \emph{attractors} play an important role in the study of their dynamics.
The maximal attractor of a dynamical system is its \emph{limit set}, or the set of points that have infinitely long chains of preimages.
Limit sets of cellular automata have been studied since the 80s, see e.g.~\cite{CuPaYu89,Ma95}.

Attractors in general, and the limit set in particular, reflect the long-term dynamics of a system.
There are related notions that correspond more strongly to the intuition of discarding transient phenomena.
For a CA $f$ and shift-invariant probability measure $\mu$ on $A^\Z$, the \emph{$\mu$-limit set} of $f$, defined by K\r{u}rka and Maass in~\cite{KuMa00}, is obtained by forbidding every word $w \in A^*$ with $\lim_{n \to \infty} f^n \mu([w]_0) = 0$.
It corresponds to observations whose probability of occurrence does not tend to 0 with time, when the initial condition is drawn from $\mu$.
The \emph{generic limit set}, defined by Milnor in~\cite{Mi85}, is the smallest closed subset of $X$ that attracts a comeager, or topologically large, subset of $X$.
Intuitively, it represents the eventual fate of all initial conditions save for a negligible set of ``pathological'' cases.
The generic limit set contains the $\mu$-limit set for sufficiently well-behaved measures $\mu$, and it is contained in the limit set.
Generic limit sets of cellular automata were first studied in~\cite{DjGu19} from the dynamical point of view, and in~\cite{To20,De21} from the computational perspective.

In this article we study the problem of determining when the generic limit set of a CA is trivial, that is, a singleton set containing a uniform configuration $\uniform{a}$ for some $a \in A$.
We call such CA \emph{generically nilpotent}.
The name comes from the classical result that a CA $f$ is nilpotent (satisfies $|f^k(A^\Z)| = 1$ for some $k \in \N$) if and only if its limit set is such a singleton set~\cite[Theorem~3.5]{CuPaYu89}.
Nilpotency is an important property of cellular automata, and several variants of it exist in the literature.
Determining whether a given CA is nilpotent is $\Sigma^0_1$-complete~\cite{Ka92}, and the analogous problem for $\mu$-limit sets, when $\mu$ is a Bernoulli measure with full support, is $\Pi^0_3$-complete~\cite[Proposition~5.6]{BoDePoSaTh15}.
We also mention \emph{asymptotic nilpotency}, which is the property that $\lim_{n \to \infty} f^n(x) = \uniform{a}$ for all $x \in A^\Z$, and \emph{unique ergodicity}, which is equivalent to $\lim_{n \to \infty} \frac{1}{n} \left| \set{0 \leq k < n}{f^k(x)_0 = a} \right| = 1$ for all $x \in A^\Z$ when $a$ is a quiescent state of $f$.
It was proved in~\cite{GuRi08}, and generalized in~\cite{Sa12,SaTo21}, that asymptotic nilpotency is equivalent to nilpotency.
On the other hand,~\cite[Proposition~48]{To15} shows that unique ergodicity is $\Pi^0_2$-complete.

We show that the class of generically nilpotent cellular automata is $\Sigma^0_2$-complete.
We also prove that even if a given CA is known to be generically nilpotent, its generic limit set cannot be determined algorithmically.
Our proofs are based on a combinatorial characterization of generic nilpotency that we give in Section~\ref{sec:Dynamics}, and the walls-and-counters technique introduced in~\cite{DePoSaTh11} for constructing cellular automata with specific asymptotic behaviors.

\section{Preliminaries}

\subsection{Definitions}

Let $A$ be a finite alphabet, whose elements are often called \emph{states}.
The one-dimensional \emph{full shift} is the set $A^\Z$ of two-way infinite sequences over $A$, called \emph{configurations}, equipped with the prodiscrete topology.
The topology is generated by the \emph{cylinder sets} $[w]_i = \set{x \in A^\Z}{x|_{[i, i+|w|)} = w}$ for $w \in A^*$ and $i \in \Z$, which are nonempty and clopen.
For a set of words $L \subset A^*$ of equal lengths, we denote $[L]_i = \bigcup_{w \in L} [w]_i$.
For $a \in A$, the all-$a$ configuration is denoted $\uniform{a}$, and such configurations are called \emph{uniform}.

The \emph{left shift} $\sigma : A^\Z \to A^\Z$ is the homeomorphism defined by $\sigma(x)_i = x_{i+1}$.
A \emph{subshift} is a topologically closed set $X \subset A^\Z$ satisfying $\sigma(X) = X$.
Alternatively, subshifts can be defined by \emph{forbidden words}: there is a set $F \subset A^*$ such that $X = A^\Z \setminus \bigcup_{w \in F} \bigcup_{i \in \Z} [w]_i$ is exactly the set configurations where none of the words in $F$ occur at any position.
For $n \in \N$, let $\lang_n(X) = \set{x|_{[0,n)}}{x \in X}$ be the set of length-$n$ words occurring in configurations of $X$.
The \emph{language} of $X$ is $\lang(X) = \bigcup_{n \in \N} \lang_n(X)$.
A standard reference for one-dimensional symbolic dynamics is~\cite{LiMa95}.

A (one-dimensional) \emph{cellular automaton}, or CA, on $A^\Z$ is a continuous function $f : A^\Z \to A^\Z$ that commutes with the shift: $f \circ \sigma = \sigma \circ f$.
By the Curtis-Hedlund-Lyndon theorem~\cite{He69}, CA are also characterized by having a finite \emph{neighborhood} $N \subset \Z$ and a \emph{local rule} $F : A^N \to A$ with $f(x)_i = F(\sigma^i(x)|_N)$ for all $x \in A^\Z$ and $i \in \Z$.
If $N \subset [-r, r]$, then $r$ is a \emph{radius} for $f$.
A state $a \in A$ is \emph{quiescent} for $f$ if $f(\uniform{a}) = \uniform{a}$.
See~\cite{Ka05} for a survey on one- and multidimensional CA.

A \emph{topological dynamical system}, or TDS, is a pair $(X, T)$ where $X$ is a compact metric space and $T : X \to X$ is continuous.
If $X \subset A^\Z$ is a subshift and $f$ is a CA on $A^\Z$, then $(X, \sigma|_X)$ and $(A^\Z, f)$ are examples of topological dynamical systems.

The \emph{limit set} of a TDS $(X, T)$ is $\limset(T) = \bigcap_{n \in \N} T^n(X)$.
The \emph{asymptotic set} $\asymset_T(x)$ of a point $x \in X$ is the set of limit points of the forward orbit $(T^n(x))_{n \in \N}$.
The \emph{realm of attraction} of $K \subset X$ is $\realm_T(K) = \set{x \in X}{\asymset_T(x) \subset K}$.
The \emph{generic limit set} $\glimset(T)$ of $(X,T)$ is the intersection of all closed subsets $K \subset X$ such that $\realm_T(K)$ is comeager in $X$.
Generic limit sets were first defined in~\cite{Mi85}.

Let $\mu$ be a Borel probability measure on $A^\Z$ and $f : A^\Z \to A^\Z$ a CA.
We can apply $f$ to $\mu$ and obtain a new probability measure $f \mu$ defined by $(f \mu)(X) = \mu(f^{-1}(X))$ for each Borel set $X \subset A^\Z$.
We say $\mu$ is \emph{shift-invariant} if $\sigma \mu = \mu$.
The $\mu$-limit set $\limset_\mu(f) \subset A^\Z$ of $f$ with respect to a shift-invariant measure $\mu$ is the subshift obtained by forbidding each word $w \in A^*$ with $\lim_{n \to \infty} f^n \mu([w]_0) = 0$.
The $\mu$-limit set was first defined in~\cite{KuMa00}.

The \emph{support} of a measure $\mu$ is the unique smallest closed set $K \subset A^\Z$ with $\mu(K) = 1$.
If $K = A^\Z$, we say $\mu$ has \emph{full support}.
The \emph{Bernoulli measure} corresponding to a distribution $p : A \to [0,1]$ with $\sum_{a \in A} p(a) = 1$ is the unique Borel probability measure $\mu_p$ on $A^\Z$ with $\mu_p([w]_i) = \prod_{j=0}^{|w|-1} p(w_j)$ for all cylinder sets $[w]_i$.
It has full support iff $p(a) > 0$ for all $a \in A$.

We say that a CA $f$ is \emph{nilpotent} if $|\limset(f)| = 1$, \emph{generically nilpotent} if $|\glimset(f)| = 1$, and \emph{$\mu$-nilpotent} for a measure $\mu$ if $|\limset_\mu(f)| = 1$.
In each case the unique configuration must be uniform.
In other words, $f$ is generically nilpotent if there exists $a \in A^\Z$ and a comeager set $U \subset A^\Z$ such that $\lim_{n \to \infty} f^n(x) = \uniform{a}$ for all $x \in U$.
It is $\mu$-nilpotent if $\lim_{n \to \infty} f^n \mu([b]_0) = 0$ for all $b \in A \setminus \{a\}$.
An alternative (and the more common) definition of nilpotency is that there exists $n \in \N$ such that $f^n(A^\Z)$ is a singleton~\cite{CuPaYu89}.


\subsection{Auxiliary results}

We list some properties of generic limit sets of cellular automata that are not evident from the definition.
The following is Proposition~4.10 from~\cite{DjGu19}.

\begin{lemma}
  \label{lem:Subshift}
  Let $f$ be a CA on $A^\Z$.
  Then $\glimset(f)$ is a nonempty $f$-invariant subshift.
\end{lemma}

The following two results are from~\cite{To20}.
The first is a combinatorial characterization of the language of the generic limit set, and the second is a tool for proving structural properties of generic limit sets.

\begin{lemma}
  \label{lem:CombChar}
  Let $f$ be a CA on $A^\Z$.
  A word $s \in A^*$ occurs in $\glimset(f)$ if and only if there exists a cylinder set $[v]_i$ such that for all $u, w \in A^*$ there exist infinitely many $t \in \N$ with $f^t([u v w]_{i - |u|}) \cap [s]_0 \neq \emptyset$.
\end{lemma}

We say that the word $v$, or the set $[v]_i$, \emph{enables} $s$ for $f$.

\begin{lemma}
  \label{lem:ForcingWords}
  Let $f$ be a CA on $A^\Z$, let $n \in \N$, and let $[v]_i \subset A^\Z$ be a cylinder set.
  Then there exists a cylinder set $[w]_j \subset [v]_i$ and $T \in \N$ such that for all $t \geq T$ we have $f^t([w]_j) \subset [\lang_n(\glimset(f))]_0$.
\end{lemma}

Words $w$ with the above property are called \emph{$\glimset(f)$-forcing}.
The lemma states that every word can be extended into a forcing one.



\section{Dynamics of generically nilpotent automata}
\label{sec:Dynamics}

Example~5.13 of~\cite{DjGu19} shows that there exist simple non-nilpotent cellular automata that are generically nilpotent.
In particular, nilpotency and generic nilpotency are not equivalent conditions.
We prove this fact for completeness, and to illustrate the use of Lemma~\ref{lem:CombChar}.
Recall that $a \in A$ is a \emph{spreading state} for a CA $f$ if the local rule $F : A^N \to A$ satisfies $F(P) = a$ whenever $P_i = a$ for some $i \in N$.

\begin{example}[Example~5.13 in~\cite{DjGu19}]
  \label{ex:min}
  Let $A = \{0,1\}$ and let $f : A^\Z \to A^\Z$ be the two-neighbor minimum CA, defined by $f(x)_i = \min(x_i, x_{i+1})$.
  The limit set $\limset(f)$ contains exactly the configurations where none of the words $1 0^n 1$ for $n \geq 0$ occur.
  In particular, $\limset(f)$ is infinite since it contains the configurations ${}^\infty 0 . 1^n 0^\infty$ for all $n \geq 0$, and $f$ is not nilpotent.
  
  Let then $s \in \lang(\glimset(f))$ be any word occurring in the generic limit set, and let $[v]_i$ enable it as per Lemma~\ref{lem:CombChar}.
  We may assume $i+|v| \geq |s|$ by extending $v$ if necessary.
  Choose $u = \epsilon$ (the empty word), $w = 0$ and $t > i+|v|$ such that $[u v w]_{i-|u|} \cap f^{-t}([s]_0)$ contains at least one configuration $x$.
  Since $x_{i+|v|} = 0$ and $0$ is a spreading state, we have $f^n(x)_j = 0$ for all $n \geq 0$ and $i+|v|-n \leq j \leq i+|v|$.
  In particular, $s = f^t(x)|_{[0,|s|)} = 0^{|s|}$.
  Since $\glimset(f)$ is a nonempty subshift, it must equal $\{ \uniform{0} \}$.
\end{example}

In~\cite{Ka92}, Kari proved that the set of nilpotent CA is $\Sigma^0_1$-complete, even within the set of CA with neighborhood $\{0,1\}$ and a spreading state.
By~\cite[Remark~1]{BoPoTh06}, every CA with a spreading state is $\mu$-nilpotent for every full-support Bernoulli measure $\mu$.
Hence nilpotent CA are $\Sigma^0_1$-complete among $\mu$-nilpotent CA.
With the same proof as Example~\ref{ex:min} we can show that every CA with neighborhood $\{0,1\}$ and a spreading state is generically nilpotent.
Hence nilpotent CA are $\Sigma^0_1$-complete among generically nilpotent CA as well.

From Proposition~6.2 and Corollary~6.5 of~\cite{DjGu19} we deduce that the $\mu$-limit set $\limset_\mu(f)$ is contained in the generic limit set $\glimset(f)$ whenever $\mu$ has full support and is $\sigma$-ergodic.
In particular, this holds when $\mu$ is a full-support Bernoulli measure, and then generic nilpotency implies $\mu$-nilpotency.
A simple example shows that the two properties are distinct.

\begin{example}[Example~5.13 in~\cite{DjGu19}]
  Let $A = \{0,1\}$ and define $f : A^\Z \to A^\Z$ by $f(x)_i = \min(x_{i+1}, x_{i+2})$.
  We claim that $f$ is $\mu$-nilpotent for every full-support Bernoulli measure $\mu$ on $A^\Z$, but not generically nilpotent.
  For the first claim, we note that $f^{-n}([1]_0) = [1^{n+1}]_n$ for all $n \geq 0$, and then $f^{-n} \mu([1]_0) = \mu([1]_0)^{n+1} \to 0$ as $n \to \infty$.
  Hence $1 \notin \lang(\limset_\mu(f))$, implying $\limset_\mu(f) = \{\uniform{0}\}$.
  Note that this is essentially the argument that a spreading state implies $\mu$-nilpotency.
  On the other hand, we have $1 \in \lang(\glimset(f))$ since the empty word enables it with $i = 0$: for all words $u, w \in A^*$ and $t \geq |w|$, there exists a configuration $x \in [u w]_{-|u|}$ with $x_j = 1$ for all $t \leq j \leq 2 t$, which implies $f^t(x)_0 = 1$.
  This also shows that a spreading state alone does not guarantee generic nilpotency.
\end{example}

In general, the limit set and the $\mu$-limit set are invariant under compositions by shifts, that is, $\limset(f \circ \sigma) = \limset(f)$ and $\limset_\mu(f \circ \sigma) = \limset_\mu(f)$ always hold.
The generic limit sets of $f$ and $f \circ \sigma$ may be distinct.
This phenomenon was studied in detail in~\cite{DjGu19}.

We now present a combinatorial characterization of generic nilpotency.

\begin{lemma}
  \label{lem:GenNil}
  Let $f$ be a CA with radius $r \geq 0$ on $A^\Z$, and let $a \in A$.
  Then $\glimset(f) = \{ \uniform{a} \}$ if and only if there exists a cylinder set $[w]_i$ and $T \in \N$ such that
  \begin{itemize}
  \item
    for all $t \geq T$ we have $f^t([w]_i) \subset [a^r]_0$, and 
  \item
    for all $n \in \N$ we have $f^{T + |A|^n}([w]_i \cap [w]_{i+r+n}) \subset [a^n]_r$.
  \end{itemize}
\end{lemma}

The first condition in particular implies that $w$ can be extended into a \emph{blocking word} for $f$, meaning a word $v \in A^*$ with some $j \in \Z$ that satisfies $\left| \set{f^t(x)|_{[0,r)}}{x \in [v]_j} \right| = 1$ for all $t \in \N$.
The intuition for the second condition is that if we ``trap'' an interval of cells between two occurrences of $w$ that force its borders to consist of $a$-states, then the deterministic dynamics of $f$ will enter a loop on this interval.
The loop must be a fixed point in which every cell contains an $a$-state, for otherwise some non-$a$ states would be found in the generic limit set.


\begin{proof}
  Suppose that $\glimset(f) = \{ \uniform{a} \}$.
  By Lemma~\ref{lem:ForcingWords} applied to the empty word $v = \epsilon$ and $n = r$, there exists $T \in \N$ and a cylinder set $[w]_i$ such that $f^t([w]_i) \subset [a^r]_0$ for all $t \geq T$.
  This is precisely the first claim.
  
  For the second claim, fix $n \in \N$ and an arbitrary configuration $x \in [w]_i \cap [w]_{i+r+n}$.
  Denote $K := \max(n+r(T+2), |i|+|w|+n+r)$ and $v = x|_{[-K,K]}$.
  Then we have $[v]_{-K} \subset [w]_i \cap [w]_{i+r+n}$.
  By the definition of $[w]_i$, this implies
  \begin{equation}
    \label{eq:borders}
    f^t([v]_{-K}) \subset [a^r]_0 \cap [a^r]_{n+r}
  \end{equation}
  for all $t \geq T$.
  Furthermore, since $r$ is a radius for $f$, we have
  \begin{equation}
    \label{eq:bottom}
    f^T(y)|_{[0, n+2r)} = f^T(x)|_{[0, n+2r)}
  \end{equation}
  for every $y \in [v]_{-K}$.

  We now prove by induction that for each $t \geq T$, the word $f^t(y)|_{[0, n+2r)} \in A^{n+2r}$ is independent of the choice of $y \in [v]_{-K}$.
  For $t = T$ this follows from~\eqref{eq:bottom}.
  Suppose then that the claim holds for some $t \geq T$.
  For each $r \leq \ell < n+r$ we have $f^{t+1}(y)_\ell = F(f^t(y)_{\ell-r}, \ldots, f^t(y)_{\ell+r})$, where $F$ is the local function of $f$.
  The right hand side is independent of $y$ by the induction hypothesis.
  For $0 \leq \ell < r$ and $n+r \leq \ell < n+2r$, we have $f^{t+1}(y)_\ell = a$ by~\eqref{eq:borders}.
  These are clearly independent of $y$.

  The above argument also shows that the sequence $s' := (f^t(x)|_{[r, n+r)})_{t \geq T}$ of length-$n$ words is eventually periodic, and the length of its pre-periodic part is less than $|A|^n$.
  Namely, each $s'(t+1)$ is determined from $s'(t)$ by the local rule $F$ and~\eqref{eq:borders}.
  The number of distinct words $s'(t) \in A^n$ for $t \geq T$ is at most $|A|^n$, and once a repetition occurs, the rest of the sequence is periodic.
  
  Again by Lemma~\ref{lem:ForcingWords}, there exists a cylinder set $C \subset [v]_{-K}$ and $T' \geq T$ with $f^t(C) \subset [a^n]_r$ for all $t \geq T'$.
  Let $y \in C$.
  Then $s'(T') = f^{T'}(y)|_{[r, n+r)} = a^n$, so that the periodic part of the sequence $s'$ is $(a^n, a^n, \ldots)$.
  Since $s'(T + |A|^n)$ is in the periodic part, we have $f^{T + |A|^n}(x)|_{[r, n+r)} = a^n$.
  The configuration $x \in [w]_i \cap [w]_{i+r+n}$ was arbitrary, so the second condition holds.

  For the converse direction, suppose that there are $[w]_i$ and $T$ that satisfy the two properties.
  Let $s \in \lang(\glimset(f))$ be arbitrary, and let $[v]_j$ enable it as per Lemma~\ref{lem:CombChar}.
  Then for all large enough $n \in \N$, there exist infinitely many $t \in \N$ with $f^t([w]_{i-n} \cap [u]_j \cap [w]_{i+n}) \cap [s]_0 = \emptyset$.
  But for each $n \geq r$ we have $f^{T + |A|^{2 n - r}}([w]_{i-n} \cap [w]_{i+n}) \subset [a^{2n-r}]_{r-n}$ by assumption, which implies $s \in a^*$.
  Because $\glimset(f)$ is a nonempty subshift of $A^\Z$, we must have $\glimset(f) = \{ \uniform{a} \}$.
\end{proof}

\section{Construction with walls and counters}
\label{sec:WallConstruction}


\subsection{Overview}

Our remaining results rely heavily on the type of construction introduced in~\cite{DePoSaTh11} and used to prove various realization results in the measure-theoretic setting~\cite{BoDeSa10,BoDePoSaTh15,HeSa18}.
We present a ``prototypical'' version of the construction in this section.

The idea is to have a designated \emph{initializer state} that can only be present in an initial configuration.
In one step, it turns into a somewhat persistent \emph{wall state}.
These walls divide a configuration into \emph{segments}.
We can achieve a high level of control over the contents of the segments by ensuring that each pair of initializers starts a process that formats the segment between them, in particular removing all walls that are not properly initialized.
When the segments are formatted, we allow them to host simulated computations, and communicate and merge with their neighbors; new walls are never created after the first time step.
In applications, we guarantee that in a typical configuration, every wall and auxiliary state will eventually disappear, so the only patterns visible in the generic limit set are those that occur inside the properly formatted segments.

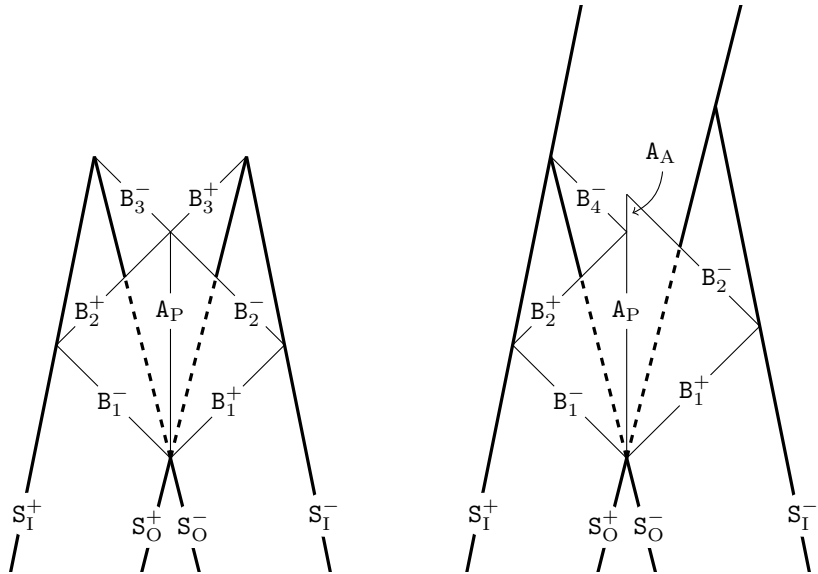
\begin{figure}[htp]
  \begin{center}
    \begin{tikzpicture}[scale=1.5]

      \begin{scope}
        \clip (-2,3) rectangle (2,8);

        \draw [very thick] (-1,0) -- (0,4);
        \draw [very thick,dashed] (0,4) -- (2/5,28/5);
        \draw [very thick] (2/5,28/5) -- (2/3,20/3);

        \draw [very thick] (1,0) -- (0,4);
        \draw [very thick,dashed] (0,4) -- (-2/5,28/5);
        \draw [very thick] (-2/5,28/5) -- (-2/3,20/3);
        \draw [very thick] (-2,0) -- ++(20/15,20/3);
        \draw [very thick] (2,0) -- ++(-20/15,20/3);

        \draw (0,4) -- ++(0,2);

        \draw (0,4) -- ++(1,1) -- ++(-5/3,5/3);
        \draw (0,4) -- ++(-1,1) -- ++(5/3,5/3);

        \node [inner sep=0cm,circle,fill=white] at (-1.25,3.5) {$\xt{S}_\mathrm{I}^+$};
        \node [inner sep=0cm,circle,fill=white] at (-0.2,3.4) {$\xt{S}_\mathrm{O}^+$};
        \node [inner sep=0cm,circle,fill=white] at (0.2,3.4) {$\xt{S}_\mathrm{O}^-$};
        \node [inner sep=0cm,circle,fill=white] at (1.35,3.5) {$\xt{S}_\mathrm{I}^-$};

        \node [inner sep=0cm,circle,fill=white] at (-0.5,4.5) {$\xt{B}_1^-$};
        \node [inner sep=0cm,circle,fill=white] at (0.5,4.5) {$\xt{B}_1^+$};
        \node [inner sep=0cm,circle,fill=white] at (-0.7,5.3) {$\xt{B}_2^+$};
        \node [inner sep=0cm,circle,fill=white] at (0.7,5.3) {$\xt{B}_2^-$};
        \node [inner sep=0cm,circle,fill=white] at (-0.3,6.3) {$\xt{B}_3^-$};
        \node [inner sep=0cm,circle,fill=white] at (0.3,6.3) {$\xt{B}_3^+$};
        
        \node [inner sep=0cm,circle,fill=white] at (0,5.3) {$\xt{A}_\mathrm{P}$};
      \end{scope}

      \begin{scope}[xshift=4cm]
        \clip (-2,3) rectangle (2,8);

        \draw [very thick] (-1,0) -- (0,4);
        \draw [very thick,dashed] (0,4) -- (7/15,28/15+4);
        \draw [very thick] (7/15,28/15+4) -- ++(3,12);
        
        \draw [very thick] (1,0) -- (0,4);
        \draw [very thick,dashed] (0,4) -- (-2/5,28/5);
        \draw [very thick] (-2/5,28/5) -- (-2/3,20/3);
        \draw [very thick] (-2,0) -- ++(30/15,30/3);
        \draw [very thick] (2+1/5,0) -- ++(-64/45,64/9);

        \draw (0,4) -- ++(0,14/6);

        \draw (0,4) -- ++(7/6,7/6) -- ++(-7/6,7/6);
        \draw (0,4) -- ++(-1,1) -- ++(1,1) -- ++(-2/3,2/3);

        \node [inner sep=0cm,circle,fill=white] at (-1.25,3.5) {$\xt{S}_\mathrm{I}^+$};
        \node [inner sep=0cm,circle,fill=white] at (-0.2,3.4) {$\xt{S}_\mathrm{O}^+$};
        \node [inner sep=0cm,circle,fill=white] at (0.2,3.4) {$\xt{S}_\mathrm{O}^-$};
        \node [inner sep=0cm,circle,fill=white] at (1.55,3.5) {$\xt{S}_\mathrm{I}^-$};

        \node [inner sep=0cm,circle,fill=white] at (-0.5,4.5) {$\xt{B}_1^-$};
        \node [inner sep=0cm,circle,fill=white] at (0.6,4.6) {$\xt{B}_1^+$};
        \node [inner sep=0cm,circle,fill=white] at (-0.7,5.3) {$\xt{B}_2^+$};
        \node [inner sep=0cm,circle,fill=white] at (0.8,5.6) {$\xt{B}_2^-$};
        \node [inner sep=0cm,circle,fill=white] at (-0.3,6.3) {$\xt{B}_4^-$};
        
        \node [inner sep=0cm,circle,fill=white] at (0,5.3) {$\xt{A}_\mathrm{P}$};
        \node [inner sep=0cm,circle,fill=white] (a) at (0.3,6.7) {$\xt{A}_\mathrm{A}$};
        \draw [->,bend left] (a) edge (0.05,6+1/6);
      \end{scope}

    \end{tikzpicture}
  \end{center}
  \caption{Comparing the values of two counters, of equal values on the left and unequal values on the right. The counter with the smaller value survives.}
  \label{fig:comparison}
\end{figure}

Figure~\ref{fig:comparison} illustrates the roles of the various signals used in the construction.
Figure~\ref{fig:init} depicts the creation of walls and formatted segments from a typical initial configuration.
The figures and the following presentation are based on~\cite[Section~3.2]{BoDePoSaTh15}, which the reader may consult for more details.
The only major difference is that in our version, wall states are created directly by initializers instead of colliding counter signals, as this simplifies some parts of the proofs.
The article~\cite{HeSa18} contains a variant in which wall states are likewise created by initializers, but the formatting process is more complex.

\subsection{The construction}

We define a radius-$1$ CA $f$ on an alphabet $A$ that realizes the prototypical segment construction.
The alphabet $A$ consists of the \emph{initializer state} $\xt{I}$, the \emph{wall state} $\xt{W}$, the \emph{empty state} $\xt{E}$, as well as states that encode \emph{signals} of various speeds.
When we say that a cell $i \in \Z$ of a configuration $x \in A^\Z$ is \emph{erased} by $f$, we mean that its next state is the empty state: $f(x)_i = \xt{E}$.
Initializers cannot be created by $f$, and every initializer will immediately turn into a wall, which is the only way a wall can be created.
A wall that originates from an initializer is called \emph{properly initialized}.

Every initializer also emits four signals: \emph{left and right inner signals} $\xt{S}^\pm_\mathrm{I}$, which move at speeds $\pm 1/5$, and \emph{left and right outer signals} $\xt{S}^\pm_\mathrm{O}$, which move at speeds $\pm 1/4$.
The outer signals can be either \emph{open} or \emph{closed}, and they start as closed.
Open signals are presesented by dashed lines in Figure~\ref{fig:comparison}.
A pair of inner and outer signals moving in the same direction, with the outer signal being the leading one, is called a \emph{counter}.
The idea is that the distance between the signals encodes the common age of the counter in unary.
A counter can appear to be older than the initial configuration, but not younger, since the signal pairs are only created by initializers and an inner signal cannot cross an outer one.

A closed outer signal erases all data it encounters, except for a closed outer signal of the opposite direction.
When two closed outer signals meet, they become open, pass through each other, and emit three new signals: \emph{left and right bouncing signals} $\xt{B}_i^\pm$ of speeds $\pm 1$, and a stationary \emph{passive anchor signal} $\xt{A}_\mathrm{P}$.
The bouncing signals come in a few different flavors ($i = 1, 2, 3, 4$), but all of them erase all data they encounter apart from inner and outer borders and the few special cases described below.
The collision of two closed outer signals first produces bouncing signals $\xt{B}_1^\pm$, which erase all data except an outer signal, which erases them, or an inner signal, which causes them to ``bounce back'' and become $\xt{B}_2^\mp$, moving in the opposite direction.

The signals $\xt{B}_2^\pm$ can pass through an open outer signal, and together with other outer signals, they are the only signals able to do so.
When this happens, the open signal becomes closed.

When two signals $\xt{B}_2^\pm$ encounter a passive anchor signal simultaneously, the latter is erased, and two new signals $\xt{B}_3^\pm$ are emitted.
When a signal $\xt{B}_3^\pm$ collides with a closed outer signal $\xt{S}^\pm_\mathrm{O}$ of the same orientation and an inner signal $\xt{S}^\mp_\mathrm{I}$ of the opposite orientation, all three are erased.
In the case that the signal $\xt{B}_3^\pm$ is produced by the collision of two equal counters, the three signals will collide simultaneously (see the left half of Figure~\ref{fig:comparison}).

When a signal $\xt{B}_2^\pm$ encounters a passive anchor signal without its symmetrical counterpart, a signal $\xt{B}_4^\mp$ is emitted in the opposite direction and the anchor becomes an \emph{active anchor signal} $\xt{A}_\mathrm{A}$.
When a signal $\xt{B}_4^\pm$ encounters a closed outer signal of the same orientation, both are erased, but an inner signal of the opposite orientation arriving at the same position is not.
If $\xt{B}_4^\pm$ is produced by the collision of unequal counters, such a three-way collision does occur (see the right half of Figure~\ref{fig:comparison}).
The only event that can erase an outer signal is a collision with $\xt{B}_3^\pm$ or $\xt{B}_4^\pm$ as described above.
When a signal $\xt{B}_2^\pm$ encounters $\xt{A}_\mathrm{A}$, both are erased.

The above scheme involves signals with fractional speeds, which can be implemented using additional states.
For example, a right inner signal with speed $1/5$ is represented by five distinct states.
They evolve cyclically and one of the transitions moves the signal by a single step.
The bouncing signals also need to remember the phases of the fractional-speed signals they bounced from, which can likewise be implemented with a finite number of additional states.

The construction ensures that when two counters collide, if they have equal values then both are destroyed, and otherwise the one with smaller value survives.
The initial configuration may contain ``rogue'' signals and walls that are not part of a counter or a comparison process, but these cannot interfere with properly initialized counters, as they will be erased by either outer signals or $\xt{B}_1^\pm$-signals.
Thus $f$ has the following property.

\begin{lemma}
  \label{lem:walls}
  Let $x \in A^\Z$ be a configuration containing at least one initializer $\xt{I}$.
  For each $i \in \Z$ with $x_i = \xt{I}$, we have $f^n(x)_i = \xt{W}$ for all $n \geq 1$.
  For each $i \in \Z$ with $x_i \neq \xt{I}$, we have $f^n(x)_i = \xt{E}$ for all $n > 5 d$, where $d = \min \set{|i-j|}{x_j = \xt{I}}$.
\end{lemma}

\begin{figure}[htp]
  \begin{center}
    \begin{tikzpicture}

      \clip (0,-1) rectangle (10,4.5);
      \draw (0,0) -- (10,0);
      
      \foreach \x/\xx in {-2/1,1/3,3/6,6/7,7/8,8/9,9/12}{
        \draw [thick,fill=black!20] (\x,0) -- (\x/2+\xx/2,\xx-\x) -- (\xx,0);
        \draw [thick] (\x,0) -- (\x/2+\xx/2,\xx*0.8-\x*0.8) -- (\xx,0);
      }
      \foreach \x in {1,3,6,7,8,9}{
        \draw [thick] (\x,0) -- (\x,4.5);
        \node [inner sep=0.05cm,circle,fill=white] at (\x,2) {$\xt{W}$};
        \node [inner sep=0.05cm,circle,fill=white] at (\x,0) {$\xt{I}$};
      }

    \end{tikzpicture}
  \end{center}
  \caption{Creation of properly initialized walls. The white regions consist of $\xt{E}$-states.}
  \label{fig:init}
\end{figure}
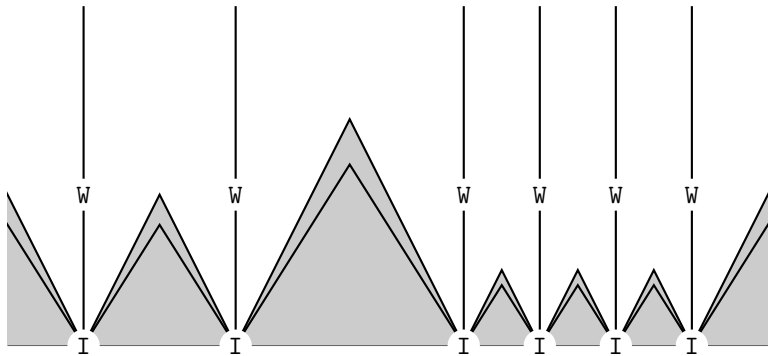

In a configuration, the space between two initializers or properly initialized walls is called a \emph{segment}.
Lemma~\ref{lem:walls} states that all segments will eventually consist of empty cells.
See Figure~\ref{fig:init}.

In our applications, we simulate computations of Turing machines inside the segments.
A simulated machine head is created to the right of each initializer (unless that cell also contains an $\xt{I}$), and the $\xt{E}$-states inside the segment are interpreted as blank tape cells.
We do not give the details of this simulation since they are standard in the literature and unimportant for our constructions, except that we must guarantee that it takes at least $5 k$ time steps for the simulated head to advance to the $k$th tape cell, for each $k \geq 1$.
Otherwise it would catch the inner signal produced by the same initializer.

\subsection{Application: finite generic limit sets}

The walls-and-counters method was used in~\cite[Theorem~6.1]{BoDePoSaTh15} to construct a CA $f$ such that for any full-support Bernoulli measure $\mu$, the $\mu$-limit set $\limset_\mu(f)$ consists of exactly two unary configurations that $f$ maps to each other.
It was noted in~\cite[Example~5.18]{DjGu19} that in this case $\glimset(f) = \limset_\mu(f)$ has the same property.
In particular, the analogue of the classical result that $\limset(f)$ is finite if and only if $f$ is nilpotent (see e.g.~\cite{Ka92}) is false for both generic limit sets and $\mu$-limit sets.

We describe the construction for completeness.
In fact, it can be easily generalized to obtain generic limit sets of any finite cardinality $1 \leq k < \infty$.
The idea is to split the empty state $\xt{E}$ into $k$ copies $\xt{E}_0, \ldots, \xt{E}_{k-1}$, and discard the wall states $\xt{W}$, so that the CA turns every initializer into $\xt{E}_1$ instead.
Each $\xt{E}_i$ becomes $\xt{E}_{i+1 \bmod k}$ under an application of the CA.
When an inner or outer signal moves to an adjacent cell or is destroyed, the cell that contained it becomes $\xt{E}_{i+1 \bmod k}$ if the neighboring cell in the opposite direction contained $\xt{E}_i$.
Then for each configuration $x$ that contains at least one initializer, for each cell $j \in \Z$ we have $f^t(x)_j = \xt{E}_{t \bmod k}$ for all large enough $t$.
It follows that $\glimset(f) = \{ \uniform{\xt{E}_0}, \ldots, \uniform{\xt{E}_{k-1}} \}$.

In~\cite[Section~4.4]{DjGu19} Djenaoui and Guillon studied dynamical systems with finite generic limit sets.
Proposition~4.16 of that article states that if $f$ is a CA with $|\glimset(f)| < \infty$, then $\glimset(f)$ consists of uniform configurations that $f$ permutes cyclically.
Hence the above construction realizes essentially all possible pairs $(\glimset(f), f|_{\glimset(f)})$, where $f$ is a CA such that $\glimset(f)$ is finite.

\section{Main results}

\begin{theorem}
  \label{thm:GenNil}
  The class of generically nilpotent cellular automata is $\Sigma^0_2$-complete.
\end{theorem}

\begin{proof}
  The combinatorial condition of Lemma~\ref{lem:GenNil} is $\Sigma^0_2$ by form.
  Hence the set of generically nilpotent CA is $\Sigma^0_2$.

  To show completeness, let $(M_k)_{k \in \N}$ be a G\"odel numbering of Turing machines, and let $P = \set{ k \in \N }{ \exists w : M_k(w) {\uparrow} }$ be the numbers of those machines that do not halt on every input.
  This set is $\Sigma^0_2$-complete.
  We many-one reduce $P$ to the set of generically nilpotent cellular automata, showing that the latter is $\Sigma^0_2$-complete as well.

  We modify the radius-$1$ cellular automaton $f : A \to A$ of Section~\ref{sec:WallConstruction}.
  The high-level idea is the following.
  Each segment will host a simulated computation of the machine $M_k$ on successive input words $w$, and the number of steps to be simulated depends on the length of the segment in such a way that $k \in P$ if and only if only finitely many lengths host a simulation that ends in a halting state.
  Whenever this happens, the segment invades a portion of its left neighbor, but in any case the eventual fate of all segments is to be filled with $\xt{E}$-states and merge with their neighbors.
  The invasions are visible in the generic limit set if and only if $k \notin P$.
  
  For each $k \in \N$, define a new alphabet $A_k = A \cup B_k$, where $B_k$ is an auxiliary alphabet that contains additional states needed to simulate a Turing machine $M'_k$ and a new state $\xt{G}$.
  We define a CA $f_k : A_k^\Z \to A_k^\Z$ such that $\glimset(f) = \{\uniform{\xt{E}}\}$ for the empty state $\xt{E} \in A$ if $k \in P$, and $\glimset(f)$ is infinite if $k \notin P$.

  The restriction $f_k|_{A^\Z}$ is exactly the walls-and-counters CA of Section~\ref{sec:WallConstruction}.
  In addition, all signals of $A$ erase all states of $B_k$ they encounter, which ensures that a configuration with infinitely many initializers will give rise to infinitely many properly initialized walls with formatted segments in between.

  Whenever $f_k$ creates a wall, it starts a simulation of the machine $M'_k$ on its right.
  This machine behaves as follows.
  First, it measures the width $n \in \N$ of the segment that hosts it.
  Then it simulates the machine $M_k$ on inputs $w_0, w_1, w_2, \ldots$, where $(w_i)_{i \in \N}$ is the enumeration of all binary words by length and then lexicographically, for a total of $n$ steps.
  If the final computation step of the simulation puts $M_k$ into a halting state, we say the segment containing $M'_k$ is \emph{bad}.
  Otherwise it is \emph{good}.
  We also define segments of length 0 to be good.

  The new state $\xt{G} \in B_k$ is by default stationary.
  When $M'_k$ has determined that its segment is good, it fills the segment with $\xt{E}$-states, except for the leftmost and rightmost cells, where it places $\xt{G}$-states.
  We extend the definition of good segments to include all segments containing $\xt{G}$-states at both ends and $\xt{E}$-states everywhere else.
  Whenever a wall has $\xt{G}$-states or walls on both sides, the central wall and any adjacent $\xt{G}$-states are immediately erased by $f_k$.
  This means that two or more neighboring good segments will merge into a larger good segment.

  If $M'_k$ determines that its segment is bad, the simulated head travels to its left end.
  It waits there until the cell at the other side of the wall contains a $\xt{G}$-state, designating the neighboring segment as good.
  Then $M'_k$ merges the two segments by erasing the wall and the $\xt{G}$-state next to it, and travels $n$ steps to the left of the position of the erased wall, or up to the left end of the combined segment, whichever is nearer, where $n$ is the length of the original bad segment.
  Note that the interior of a good segment contains only $\xt{E}$-states, so the head can travel through it undisturbed.
  Then the head travels to the right end of the combined segment, erasing all data from the combined segment as it goes.
  It places a $\xt{G}$-state next to the right wall and erases itself.
  This means that a bad segment will merge with a good segment to its left, and the combined segment then becomes good.
  In the process a simulated head of $M'_k$ travels $n$ cells deep into the good segment.
  
  This concludes the definition of $f_k$.
  We note that $f_k$ can be realized with radius $r = 2$.

  Suppose that $k \in P$.
  Since there is an input on which $M_k$ never halts, the set of lengths that result in bad segments is finite, so it has an upper bound $N \in \N$.
  We claim that the word $w = \xt{I} \xt{E}^{N+3} \xt{I} \in B_k^{N+5}$ satisfies the conditions of Lemma~\ref{lem:GenNil} with $i = -2$ and some $T \in \N$.
  Consider a configuration $x \in [w]_{-2}$.
  The two $\xt{I}$-states $x_{-2}$ and $x_{N+2}$ delimit a segment, which is properly formatted in $f_k^{5 N}(x)$.
  After the simulated machine $M'_k$ in this segment has finished its computation, it becomes good, since it is strictly longer than $N$ cells.
  Let $T$ be the time step at which this happens.
  Then $f^T(x)|_{[0,1]} = \xt{E}^2$, and we claim that $f^t(x)|_{[0,1]} = \xt{E}^2$ holds for all $t \geq T$.

  For $t \geq T$, let $S_t \subset \Z$ be the segment of $f^t(x)$ containing the coordinate $0$.
  If $S_t$ merges with one or more good segments, the resulting segment $S_{t+1}$ is also good, and the merging process does not involve the interior cells of $S_t$.
  If $S_t$ merges with a bad segment to its right, then a simulated head of $M'_k$ will travel some $n \leq N$ cells into $S_t$.
  Since the rightmost coordinate of $S_t$ is initially $N+1$ and cannot decrease, the head will not reach the coordinate $1$.
  This shows $f^t(x)|_{[0,1]} = \xt{E}^2$ for all $t \geq T$.

  For the second condition of Lemma~\ref{lem:GenNil}, take any $n \in \N$ and a configuration $x \in [w]_{-2} \cap [w]_{n-2}$.
  The word $x|_{[-2, n+N+3]}$ contains some number of $\xt{I}$-states that divide it into segments.
  Each of those segments will eventually be properly formatted and become either good or bad.
  The first and last segments will become good.
  It follows that the segment of $f^t(x)$ containing the coordinate $0$ will always be good, and it will merge with the segment to its right as soon as the latter becomes good or bad.
  Hence, for large enough $t$ we will have $f^t(x)|_{[0, n+N+1]} = \xt{E}^{n+N+2}$, which implies the second condition of Lemma~\ref{lem:GenNil}.
  Thus $f_k$ is generically nilpotent.
  
  Suppose now that $k \notin P$.
  Then there exist infinitely many lengths that result in bad segments.
  For each $n \geq 0$, denote $L_n = \set{\xt{E}^n b}{b \in B_k \setminus \{\xt{E}\}} \subset B_k^{n+1}$.
  We show that $L_n \cap \lang(\glimset(f)) \neq \emptyset$ for all $n \geq 0$.
  More strongly, we show that some word in $L_n$ is enabled by the empty word, using the terminology of Lemma~\ref{lem:CombChar}, which is equivalent to the condition that for every cylinder set $C$ there are infinitely many $t \in \N$ with $f^t(C) \cap [L_n]_0 \neq \emptyset$.
  This implies that $\glimset(f_k)$ is an infinite subshift, so $f_k$ is not generically nilpotent.

  Let $[w]_i$ be a cylinder set.
  We may assume $i \leq 0$ and $|w| \geq |i|+n$ by extending $w$ if necessary.
  Choose a number $N \geq |w|-n$ that results in a bad segment and consider the configuration
  $x = {}^\infty \xt{E} \xt{I} \xt{I} . w \xt{I} \xt{E}^N \xt{I} \xt{E}^\infty$,
  where the dot denotes coordinate $i$, so that $x \in [w]_i$.
  We claim that if $N$ is large enough, then $f_k^t(x) \in [L_n]_0$ for some $t > N$.
  As we have infinitely many choices for $N$, this implies that $f^t([w]_i) \cap [L_n]_0 \neq \emptyset$ for infinitely many $t$.

  Every finite segment in $x$ will eventually become either good or bad.
  The leftmost segment has length 0, and we have guaranteed that such segments are good.
  Therefore it will merge with its right neighbor once that segment has determined its goodness status, and the combined segment will become good.
  In this way all segments in $x$ will eventually merge into one.
  
  Suppose $N$ is so large that the rightmost segment becomes bad only after all the other segments have merged into a good segment.
  This requires more than $N$ applications of $f_k$, since it takes at least $N$ steps for the simulates head of $M'_k$ to measure the length of the segment.
  Since we assumed $N \geq |w|-n$, the head of $M'_k$ will travel to the left, and at some time step $t > N$ enter the coordinate $n$.
  Then $f_k^t(x)|_{[0,n]} \in L_n$, which is what we claimed.
  This shows that $\glimset(f_k)$ is infinite.
\end{proof}

In~\cite{BoDePoSaTh15}, the authors consider the complexity of $\mu$-nilpotent cellular automata also within specific subclasses.
They prove that for a $\sigma$-ergodic measure $\mu$ with full support, $\mu$-nilpotent CA are $\Pi^0_1$ within the class of CA with a \emph{persistent state} (a state $a \in A$ with $f(x)_i = a$ whenever $x_i = a$), and $\Sigma^0_2$ within the class of CA with an equicontinuity point.
The latter result is somewhat analogous to Theorem~\ref{thm:GenNil}: as discussed after Lemma~\ref{lem:GenNil}, every generically nilpotent CA admits a blocking word, and hence an equicontinuity point~\cite{Ku97}.

Delacourt proved in~\cite{De21} that every nontrivial property of generic limit sets of cellular automata is undecidable.
More strongly, they are $\Pi^0_1$-hard or $\Sigma^0_1$-hard.
It was left open whether this complexity bound is tight, that is, whether there exist $\Pi^0_1$ properties of generic limit sets.
We do not resolve this problem, but Theorem~\ref{thm:GenNil} shows that some properties are $\Sigma^0_2$-complete.

If $f : A^\Z \to A^\Z$ is a nilpotent CA, then $\limset(f) = \{\uniform{a}\}$ for some $a \in A$.
Determining the state $a$ from the local rule of $f$ is simple: it is the sole quiescent state of $f$.
We show that it is impossible to algorithmically determine the generic limit set of a given generically nilpotent CA.

\begin{theorem}
  \label{thm:DetermineState}
  Given a finite alphabet $A$, distinct states $a, b \in A$, and a cellular automaton $f : A^\Z \to A^\Z$ that satisfies either $\glimset(f) = \{\uniform{a}\}$ or $\glimset(f) = \{\uniform{b}\}$, it is undecidable whether the first case holds.
\end{theorem}

\begin{proof}
  We reduce the halting problem of Turing machines on empty input to the problem in the theorem statement.
  Let $M$ be a turing machine.
  We define a cellular automaton $f_M$ as follows.
  As in the proof of Theorem~\ref{thm:GenNil}, we modify the radius-$1$ CA $f : A^\Z \to A^\Z$ of Section~\ref{sec:WallConstruction} by extending the alphabet into $A_M := A \cup B_M$.
  The auxiliary alphabet $B_M$ contains new states that are used to simulate a Turing machine $M'$, as well as new states $\xt{G}$ and $\xt{H}$.
  The state $\xt{G}$  plays a similar role as in the earlier proof, marking segments that are ready to merge with their neighbors.
  The state $\xt{H}$ spreads over $\xt{E}$-states in both directions, but not over other states of $A_M$.
  
  As a wall is created by an initializer, a simulated computation of $M'$ is started on its right.
  This machine measures the length $n$ of the segment containing it and simulates $M$ on empty input for $n$ steps.
  Then it writes $\xt{G}$-states on the leftmost and rightmost cells of the segment.
  If $M$ halted during the simulated computation, then $M'$ writes an $\xt{H}$-state somewhere in the interior of the segment.
  Then, regardless of whether $M$ halted, $M'$ erases itself.
  As before, if a wall is surrounded by walls or $\xt{G}$-states, it erases itself and any adjacent $\xt{G}$-states.
  
  This concludes the definition of $f_M$.
  It can be implemented with radius $r = 2$.

  We claim that $\glimset(f_M) = \{\uniform{\xt{H}}\}$ if $M$ halts on empty input, and $\glimset(f_M) = \{\uniform{\xt{E}}\}$ if it does not.
  Suppose first that $M$ halts after some $N$ steps.
  We claim that the word $w = \xt{I} \xt{E}^{N+4} \xt{I}$, $i = -2$ and $a = \xt{H}$ satisfy the two conditions of Lemma~\ref{lem:GenNil}, showing that $\glimset(f_M) = \{\uniform{\xt{H}}\}$.
  
  For the first condition, let $x \in [w]_{-2}$ be arbitrary.
  The initializers at $x_{-2}$ and $x_{N+3}$ delimit a segment of length $N+4$.
  In this segment, the machine $M'$ is simulated, which in turn simulates $M$ until it halts.
  Then it writes $\xt{G}$-states at coordinates $-1$ and $N+2$ and an $\xt{H}$ somewhere in between before erasing itself.
  The $\xt{H}$ will spread over the cells at coordinates $0, 1, \ldots, N+1$.
  After this, their contents will never change, as the merging of segments only affects the cells adjacent to the walls.
  In particular, we have $f_M^t(x)|_{[0,1]} = \xt{H}^2$ for all large enough $t$, and thus the first condition holds.

  For the second condition, take any $n \in \N$ and a configuration $x \in [w]_{-2} \cap [w]_{n-2}$.
  As in the proof of Theorem~\ref{thm:GenNil}, the word $x|_{[-2, n+N+4]}$ contains some number of $\xt{I}$-states that divide it into segments, which will eventually merge into one large segment.
  This segment contains at least one $\xt{H}$, since the original leftmost segment has length $N+4$.
  Hence $f_M^t(x)|_{[-2, n+N+4]} = \xt{H}^{n+N+6}$ for all large enough $t$, and the second condition holds.

  Suppose then that $M$ does not halt on empty input.
  We show that $w = \xt{I} \xt{E}^4 \xt{I}$, $i = -2$ and $a = \xt{E}$ satisfy the two conditions of Lemma~\ref{lem:GenNil}, implying $\glimset(f_M) = \{\uniform{\xt{E}}\}$.
  Let $x \in [w]_{-2}$ be arbitrary.
  The segment delimited by $x_{-2} = x_3 = \xt{I}$ will host a simulation of $M'$, which simulates $M$ for 4 steps and writes $\xt{G} \xt{E} \xt{E} \xt{G}$ on its cells.
  In general, every properly formatted length-$n$ segment in $x$ will eventually be filled with $\xt{E}$-states, except for the two bordermost cells.
  Merging two such segments preserves this property, and the walls-and-counters construction ensures that segments can only be modified by merging.
  Thus $f_M^t(x)|_{[0,1]} = \xt{E}^2$ for all large enough $t$, and the first condition holds.

  The proof of the second condition is analogous to the case of $M$ halting, except that no $\xt{H}$-symbols are produced in any segments.
\end{proof}

\section{Future work}

In the proof of Theorem~\ref{thm:DetermineState}, we actually show that the problem of determining the generic limit set of a generically nilpotent CA is $\Sigma^0_1$-hard.
By symmetry, it is also $\Pi^0_1$-hard.
This implies that it is not complete for either class.
It is $\Sigma^0_2$ by Lemma~\ref{lem:GenNil}, but not $\Sigma^0_2$-complete for the same reason as above.
We leave its exact complexity open.

In this article we have only considered one-dimensional cellular automata.
The proof of Lemma~\ref{lem:GenNil} cannot be directly generalized to the multidimensional case, since it relies on the ability of blocking words to cut a one-dimensional configuration into two independently evolving halves.
Thus, we do not know if our results hold for two- and higher-dimensional CA.
A two-dimensional version of the walls-and-counters CA of Section~\ref{sec:WallConstruction} was presented in~\cite{DeHe15}.

Finally, we have concentrated on generic limit sets consisting of a single configuration.
Our results probably apply to CA with finite generic limit sets with mostly the same proofs, save for additional technical details.

\bibliographystyle{plain}
\bibliography{GenericsBib}

\end{document}